%% file: doubleint.tex
\documentclass[a4paper,twoside,12pt]{article}
\usepackage{amsmath,amsfonts,amsthm,amsbsy,amssymb,amscd}
\usepackage[all]{xy}
\newcommand{\ZZ}{\mathbb{Z}} 
\newcommand{\PP}{\mathbb{P}}
\newcommand{\QQ}{\mathbb{Q}}

\newcommand{\FQ}{\mathbb{F}_{q}}
\newcommand{\CC}{\textbf{C}}

\newcommand{\OO}{\mathcal{O}}
\newcommand{\HH}{\mathcal{H}}
\newcommand{\EE}{\mathcal{E}}
\newcommand{\FF}{\mathcal{F}}
\newcommand{\VV}{\mathcal{V}}
\newcommand{\II}{\mathcal{I}}

\newcommand{\Tp}{\mathcal{T}_\mathfrak{p}}
\newcommand{\Tinf}{\mathcal{T}_\infty}
\newcommand{\Tree}{\mathcal{T}}
\newcommand{\KK}{\mathcal{K}}
\newcommand{\Ac}{\mathcal{A}}
\newcommand{\Cc}{\mathcal{C}}
\newcommand{\LL}{\mathcal{L}}

\newcommand{\nn}{\mathfrak{n}}
\newcommand{\fp}{\mathfrak{p}}
\newcommand{\fq}{\mathfrak{q}}

\newcommand{\ud}{\underline}
\newcommand{\rsp}{{\rm sp}}

\newcommand{\co}{{\underline{\rm{H}}}}
\newcommand{\cuco}{{\underline{\rm H}}_{!}}
\newcommand{\hc}[3]{{\underline{\rm H}}_{!} (\mathcal{T}_{#1},#2)^{#3}}

\newcommand{\gl}{\mbox{GL}}
\newcommand{\pgl}{\mbox{PGL}}
\newcommand{\SL}{\mbox{SL}}

\newcommand{\GaB}{{\Gamma_0^\mathfrak{p}(\nn\infty)}}

\newcommand{\ord}{\mbox{ord}}

\newcommand{\Gamtil}{{\tilde{\Gamma}}}
\newcommand{\teit}{\mu_{{\rm Teit}}}
\newcommand{\mint}{\times\!\!\!\!\!\!\int}

\newcommand{\cxy}{c\{\infty\to 0\}}
\newcommand{\mcxy}{\mu_c\{\infty\to 0\}}

\author{{\it Hilmar Hauer} (Nottingham) and {\it Ignazio Longhi} (Padova)}
\title{Teitelbaum's exceptional zero conjecture in the function field case}

\begin{document}
\maketitle
\tableofcontents 
\pagestyle{headings} 
\theoremstyle{plain}

\newtheorem{Def}{Definition}[section]
\newtheorem{Prop}[Def]{Proposition}
\newtheorem{Thm}[Def]{Theorem} 
\newtheorem{Lemma}[Def]{Lemma}
\newtheorem{Cor}[Def]{Corollary}
\newtheorem{Conj}[Def]{Conjecture}
\newtheorem{Ex}[Def]{Example} 
\newtheorem*{Rem}{Remark}

\include{intro}
\include{notation}
\include{prelim}
\include{integr}

\include{Da_conj}

\include{bibl}
\end{document}

%% file: intro.tex
\section*{Introduction}

Let $E$ be an elliptic curve over $\QQ$ of conductor $N$ and $p$
a prime number. 
Since $E$ is modular, it corresponds to a cusp form of weight 2, 
i.e.\ an analytic function $f(z)$ on the complex upper half plane 
$\mathcal{H}$ s.t.\ 
$$
f\left(\frac{az+b}{cz+d}\right)=(cz+d)^2f(z)
\textrm{ for all }
\begin{pmatrix} a & b \\ c & d \end{pmatrix} \in \Gamma_0(N),
$$
and it fulfills suitable growth conditions at the points in
$\PP^1(\QQ)$ (the so-called cusps).
Mazur, Tate and Teitelbaum construct in \cite{mtt}
a $p$-adic $L$-function $L_p(E,s)$ using a $p$-adic measure
associated to $f(z)$ via the modular symbols
$$
\Psi_E:{\rm Div}^0(\PP^1(\QQ))\longrightarrow\QQ
$$
associated to $E$.

Furthermore, they propose $p$-adic versions
of the Birch and Swinnerton-Dyer conjecture.
It turns out that $L_p(E,s)$ always vanishes at the central
point $1$ if $E$ has split multiplicative reduction at
$p$. In this case they conjecture (the so-called 
\textit{exceptional zero conjecture}):
$$
L'_p(E,1)=\frac{\log_p(q_E)}{\ord_p(q_E)}L_{\rm norm}(E,1),
$$
where $q_E$ is the Tate period of $E$ at $p$, $\log_p$ the $p$-adic
logarithm with $\log_p(p)=1$ and $L_{\rm norm}(E,s)$ a suitable 
normalisation of the complex $L$-function.
This was proved by Greenberg and Stevens \cite{gs}.
% and independently by \textbf{KKT??????}. 
The proof is very technical and makes use of Hida's theory of
$p$-adic $L$-functions.

In their paper \cite{mt}, Mazur and Tate state slightly stronger,
refined versions of their conjectures, avoiding the construction of
a $p$-adic $L$-function.
In particular, they give an exponentiated version of the exceptional
zero conjecture:
$$
\prod_{\stackrel{a\textrm{ mod }p^n}{(a,p)=1}}
a^{{\rm ord}_p(q_E)\cdot \Psi_E((a/p^n)-(i\infty))}
\equiv
\left(\frac{q_E}{p^{{\rm ord}_p(q_E)}}\right)^{\Psi_E((0)-(i\infty))}
\textrm{ in }(\ZZ/p^n\ZZ)^*
$$
for every $n\geq 1$ (proved by de Shalit \cite{ds} under certain
assumptions).

\medskip
Now let $E$ be an elliptic curve over $F=\FQ(T)$ of conductor
$\nn\infty$ with split multiplicative reduction at the place
$\infty=\frac{1}{T}$ and $F_\infty$ the completion of $F$ at $\infty$.
By Drinfeld's work \cite{dr},
%(combined with results of Weil, Grothendieck,
%Jacquet-Langlands, Deligne and Zarhin)
$E$ has a uniformisation
$$
\overline{M}_0(\nn)\longrightarrow E
$$
by a Drinfeld modular curve and it corresponds to a
Drinfeld cusp form on 
$\Omega_\infty=\PP^1(\CC_\infty)-\PP^1(F_\infty)$
($\CC_\infty$ a completion of an algebraic closure of $F_\infty$) 
or, equivalently,
to a cuspidal harmonic cochain $c_\infty$ on the Bruhat-Tits tree
$\Tinf$ (with values in $\QQ$).
Teitelbaum (\cite{t3}) defines the group of modular symbols
to be
$$
\mathcal{M}:={\rm Div}^0(\PP^1(F))
$$
and denotes the divisor $(r)-(s)$ by $[r,s]$.
He defines a map
$$
\begin{array}{ccc}
\mathcal{M} & \longrightarrow & \QQ \\
{[r,s]} & \longmapsto & [r,s]\cdot c_\infty,
\end{array}
$$
where $[r,s]\cdot c_\infty$ is the sum of the values of $c_\infty$
along the axis connecting the ends corresponding to $r$ and $s$ on
$\Tinf$.

If $E$ has split multiplicative reduction
at a further place $\mathfrak{p}$ different from $\infty$, he defines 
a measure on $\OO_\mathfrak{p}$ analogous to the classical case.
(Here $\OO_\mathfrak{p}$ is the ring of integers of the completion of $F$ at
$\mathfrak{p}$.)
However, the absence of a logarithm obstructs the definition
of a $\mathfrak{p}$-adic $L$-function analogous to the number field
case. But it is still possible to formulate the refined
exceptional zero conjecture:
$$
\left( \frac{q_\mathfrak{p}}{\pi_\mathfrak{p}^{\textrm{ord}_\mathfrak{p}(q_\mathfrak{p})}} 
\right)^{[0,\infty]\cdot c_\infty}
=
\left( \lim_{n\to\infty} 
\prod_{\stackrel{a\textrm{ mod }\pi_\mathfrak{p}^n}{(a,\pi_\mathfrak{p})=1}}
a^{[a/\pi_\mathfrak{p}^n,\infty]\cdot c_\infty} 
\right)^{\textrm{ord}_\mathfrak{p}(q_\mathfrak{p})},
$$
where $q_\mathfrak{p}$ is the Tate period associated to $E$ at $\mathfrak{p}$.
This is carried out in \cite{t3}, where also some numerical
evidence can be found.

\bigskip
In this paper we prove Teitelbaum's conjecture up to a root of unity.
Therefore we can regard
our result as an analogue of the Greenberg-Stevens formula.
We achieve this by developing an analogue of Darmon's integration theory
on $\mathcal{H}_p\times\mathcal{H}$, where 
$\mathcal{H}_p=\PP^1(\CC_p)-\PP^1(\QQ_p)$ is the $p$-adic upper half plane
\cite{da}. This is done in sections 2.1 and 2.2 for arbitrary global function
fields and choices of $\fp$ and $\infty$, i.e., in a more general setting
than in \cite{t3}. The function field situation 
is considerably easier than the number field case. Since both places
are non-archimedean, we can work with certain harmonic cochains
on the product of the Bruhat-Tits trees at $\mathfrak{p}$ and $\infty$:
$$
\Tree=\Tp\times\Tinf.
$$
By this we mean the set product of the vertices/edges. The symmetry of this
situation will prove very useful for later calculations.

Since $E$ has split multiplicative reduction at $\mathfrak{p}$ and $\infty$,
it corresponds to two harmonic cochains $c_\mathfrak{p}$ and $c_\infty$ on
$\Tp$ and $\Tinf$, respectively. In section 2.3 we compare them by
relating them to a certain space of automorphic newforms.

Darmon defines a period $I_\psi\in\CC_p^*$ (for details compare
\cite{da}). The order at $p$
and the $p$-adic logarithm of this period are closely related to special
values of a certain partial $L$-function and the first derivative of
the $p$-adic $L$-function attached to $E/\QQ$, respectively.
He obtains the following reformulation of the result of Greenberg
and Stevens:
$$
\log_p(I_\psi)=\frac{\log_p(q_E)}{\ord_p(q_E)}\ord_p(I_\psi).
$$
We define an analogous period in section 3.1. Furthermore, we show an
exponentiated version of the above formula (theorem \ref{da_thm}):
$$
I_\psi=\zeta\cdot q_\mathfrak{p}^{\frac{\nu_\mathfrak{p}(I_\psi)}{\nu_\mathfrak{p}(q_\mathfrak{p})}},
$$
where $\zeta$ is a root of unity.

In the final section we assume the setting in \cite{t3} and
show that Teitelbaum's conjecture
is equivalent (up to a root of unity) to the above theorem.

\bigskip
\textit{Acknowledgements.}
We are indebted to Michael Spie{\ss} and Henri Darmon for
suggesting this problem and many helpful remarks. Furthermore,
Longhi would also like to thank 
Gebhard B\"ockle, Ernst-Ulrich Gekeler, 
Matteo Longo, Richard Pink and Hans-Georg R\"uck for useful
discussions and remarks.
While working on this paper,
Hauer was partially supported by EPSRC grant GR/M89560.
Longhi was supported by post-doc scholarships of the Arithmetic
Algebraic Geometry and the Galois Theory and Explicit Methods in
Arithmetic networks, and by a post-doc scholarship of the 
Universit\`a di Padova.

%% file: notation.tex
\section*{Notations}

$F$ is the function field of a geometrically connected 
smooth projective algebraic curve
$\mathcal{C}$ over the finite field $\FQ$ of characteristic $p$.
For a closed point $\mathfrak{p}\in\mathcal{C}$ denote by
\begin{itemize}
\item 
$A_\mathfrak{p}=\Gamma(\mathcal{C}-\{\mathfrak{p}\},\OO_{\mathcal{C}})$
the ring of regular functions on $\mathcal{C}-\{\mathfrak{p}\}$
and more generally,
$A_S=\Gamma(\mathcal{C}-S,\OO_{\mathcal{C}})$ the ring of regular functions
on $\mathcal{C}-S$ for a finite set $S$ of closed points,
\item 
$\nu_\mathfrak{p}$ the corresponding valuation of $F$ at $\fp$,
\item
$\pi_\mathfrak{p}$ a uniformiser,
\item
$F_\mathfrak{p}$ the completion of $F$ at $\mathfrak{p}$ with valuation ring
$\OO_\mathfrak{p}$ and maximal ideal $\mathfrak{m}_\mathfrak{p}$,
\item
$\CC_\mathfrak{p}$ the completion of an algebraic closure of $F_\mathfrak{p}$,
\item
$\Omega_\mathfrak{p}=\PP^1(\CC_\mathfrak{p})-\PP^1(F_\mathfrak{p})$ Drinfeld's upper
half plane and $\Omega_\mathfrak{p}^*=\Omega_\mathfrak{p}\cup\PP^1(F)$.
\item
$\Ac$ is the adele ring of $F$ with ring of integers $\OO$. These
decompose into a $\mathfrak{p}$-part and a "finite" part with respect to $\mathfrak{p}$:
$$
\Ac=\Ac_{f,\mathfrak{p}}\times F_\mathfrak{p}
\textrm{\ \ \ and\ \ \ }
\OO=\OO_{f,\mathfrak{p}}\times\OO_\mathfrak{p}.
$$
Elements of $\Ac$ are denoted by $\underline{x}=(x_\mathfrak{q})_\mathfrak{q}$.
\end{itemize}

\noindent
A subgroup $\Gamma$ of $\gl_2(F)$ is \textit{arithmetic} 
(w.r.t.\ $\mathfrak{p}$) if it
is commensurable with $\gl_2(A_\mathfrak{p})$. 
For such a group, we write
$$
\overline{\overline{\Gamma}}:=\Gamma/\Gamma\cap {\rm Center}(\gl_2(F))
$$
and
$$
\overline{\Gamma}:=\Gamma^{{\rm ab}}/{\rm torsion}.
$$
Sometimes we use $G$ and $Z$ to denote $\gl_2$ and its centre. 

\noindent
For an ideal $\nn$ of $A_\mathfrak{p}$, we define the
\textit{Hecke congruence subgroup associated to}
$\nn$ (w.r.t.\ $\mathfrak{p}$): 
$$
\Gamma_0^\mathfrak{p}(\nn):=
\left\{
\gamma\in\mbox{GL}_2(A_\mathfrak{p}) : \gamma\equiv 
\left(\begin{array}{cc}\ast & \ast \\ 0 & \ast \end{array}\right)
\textrm{ mod } \nn
\right\}.
$$

%% file: prelim.tex
\section{Preliminaries}

\subsection{The Bruhat-Tits Tree for $\pgl_2(F_\mathfrak{p})$}

The vertices of the tree $\Tp$ are defined to be homothety classes
of $\OO_\mathfrak{p}$-lattices (i.e.\ free $\OO_\mathfrak{p}$-modules of rank $2$) in
$F_\mathfrak{p}^2$. Two such classes $[L],[L']$ are connected by an edge if
there exists $L''\in [L']$ such that $L''\subseteq L$ and $L/L''$ has
length $1$ as $\OO_\mathfrak{p}$-module, i.e.\ 
$L/L''\cong\OO_\mathfrak{p}/\pi_\mathfrak{p}\OO_\mathfrak{p}$.

We denote by $\EE(\Tp),\vec{\EE}(\Tp),\VV(\Tp)$ respectively the set
of edges, oriented edges and vertices of $\Tp$.
We consider the following action of $\gl_2(F_\mathfrak{p})$
from the left on $\Tp$:
$$
\gamma_*[L]:=[L\gamma^{-1}].
$$
We note that this action is different from the one considered in \cite{se}.
It induces the following identifications:
$$
\begin{array}{ccc}
\gl_2(F_\mathfrak{p})/\gl_2(\OO_\mathfrak{p})\cdot {\rm Z}(F_\mathfrak{p})
    & \tilde{\longrightarrow} & \VV(\Tp) \\
\gamma & \longmapsto & \gamma_* v_0
\end{array}
$$
and
$$
\begin{array}{ccc}
\gl_2(F_\mathfrak{p})/\II_\mathfrak{p}\cdot {\rm Z}(F_\mathfrak{p})
    & \tilde{\longrightarrow} & \vec{\EE}(\Tp) \\
\gamma & \longmapsto & \gamma_* e_0,
\end{array}
$$
where ${\rm Z}$ is the centre of $\gl_2$,
$$
\II_\mathfrak{p}=\left\{\begin{pmatrix} a & b \\ c & d \end{pmatrix}\in\gl_2(\OO_\mathfrak{p}):
c\equiv 0\textrm{ mod }\mathfrak{p}\right\}
$$
the Iwahori group,
$v_0=[\OO_\mathfrak{p}^2]$, $v_{-1}=[\pi_\fp\OO_\mathfrak{p}\oplus\OO_\mathfrak{p}]$,
$e_0$ the edge with origin $o(e_0)=v_{-1}$ and terminus $t(e_0)=v_0$,
and $\overline{e}_0$ is the edge opposite to $e_0$.
We also note that $\gl_2(F_\fp)$ acts canonically from the left on
$\gl_2(F_\mathfrak{p})/\gl_2(\OO_\mathfrak{p})\cdot {\rm Z}(F_\mathfrak{p})$
and
$\gl_2(F_\mathfrak{p})/\II_\mathfrak{p}\cdot {\rm Z}(F_\mathfrak{p}).$

\begin{Lemma} \label{Tp}
Each edge $e$ of $\Tp$ can uniquely be written as a product $\gamma e_0$,
where $\gamma$ is an element of the following disjoint union of sets:
$$
\begin{array}{l}
\left\{ \begin{pmatrix} \pi_\mathfrak{p}^n & u \\ 0 & 1 \end{pmatrix} : 
        n\in\ZZ \textrm{, } u\in F_\mathfrak{p}/\pi_\mathfrak{p}^n\OO_\mathfrak{p} \right\} \\
\dot{\cup}\textrm{\ }
\left\{ \begin{pmatrix} \pi_\mathfrak{p}^n & u \\ 0 & 1 \end{pmatrix} 
        \begin{pmatrix} 0 & 1 \\ \pi_\mathfrak{p} & 0 \end{pmatrix} : 
        n\in\ZZ \textrm{, } u\in F_\mathfrak{p}/\pi_\mathfrak{p}^n\OO_\mathfrak{p} \right\}.
\end{array}
$$
\end{Lemma}

\begin{proof}
This follows easily from the above identification (see \cite{gek1}).
\end{proof}

The ends of $\Tp$ (=infinite half-lines, equivalent if coinciding up to
a finite number of edges) correspond to 
$\partial\Omega_\mathfrak{p}=\PP^1(F_\mathfrak{p})$. We label the path 
$A(\infty,0)$ from $\infty$ to $0$ by
$(v_i)_{-\infty<i<\infty}$ resp.\ $(e_i)_{i\in\ZZ}$,
where $t(e_i)=v_i$.
In the following, we will call $v_0$ the \textit{base vertex}, 
$e_0$ the \textit{base edge}
and $A(\infty,0)$ the \textit{base axis} of $\Tp$.
In particular, $e_i$ (and $v_i$) is represented by
$\begin{pmatrix} \pi_\fp^i & 0 \\ 0 & 1 \end{pmatrix}$ in lemma \ref{Tp}. 
(For more detail see \cite[1.3]{gr} and \cite{se}.)

\bigskip
Let $M$ be an abelian group. A harmonic cochain with values in $M$ is
a map $\varphi:\vec{\EE}(\Tp)\rightarrow M$ that is alternating, i.e.\ 
$\varphi(\overline{e})=-\varphi(e)$ for all $e\in\vec{\EE}(\Tp)$ and harmonic, i.e.\
$$
\sum_{t(e)=v}\varphi(e)=0
$$ 
for all $v\in\VV(\Tp)$.
We denote the group of $M$-valued harmonic cochains by $\co(\Tp,M)$.

The identification of $\partial\Omega_\mathfrak{p}$ with the ends
of $\Tp$ induces a bijection between $\co(\Tp,M)$ and 
$\textrm{Meas}_0(\partial\Omega_\mathfrak{p},M)$, the set of
$M$-valued measures (i.e., finitely additive functions on
compact open subsets) of total mass $0$ on $\partial\Omega_\mathfrak{p}$.
If $\varphi$ is a harmonic cochain, then
$$
\mu(U(e)):=-\varphi(e)
$$
defines a measure of total mass $0$, where
the set $U(e)$ consists of all ends going through $e\in\vec{\EE}(\Tp)$.

\medskip
Let $\Gamma$ be an arithmetic subgroup of $\gl_2(F)$. Then the group
of $M$-valued cuspidal harmonic cochains for $\Gamma$ is the following:
$$
\cuco(\Tp,M)^\Gamma:=\left\{\varphi\in\co(\Tp,M):
\begin{array}{l}
\varphi \textrm{ is }\Gamma \textrm{-invariant and} \\
\textrm{of finite support modulo } \Gamma
\end{array}
\right\}.
$$
The group $\gl_2(F_\mathfrak{p})$ acts on $\cuco(\Tp,M)^\Gamma$ by
$$
\gamma\ast\varphi:=\varphi\circ\gamma^{-1}.
$$
This action is compatible with its action on $\Tp$.

There is an injection
$$
j:\overline{\Gamma}\longrightarrow\cuco(\Tp,\ZZ)^\Gamma,
$$
with finite cokernel,
such that for a chosen vertex $v\in\VV(\Tp)$, $j(\alpha)(e)$ 
counts the number of $\gamma e$ lying in $A(v,\alpha v)$ as
$\gamma$ varies in $\overline{\overline{\Gamma}}$ (\cite[3.3.3]{gr}).
It is proved to be bijective if $A_\fp$ is a polynomial ring and
$\Gamma$ is a Hecke congruence group (\cite[thm.\ 1.9]{gek}).
Although unknown to be surjective in general, no case is known
where this map fails to be bijective.

If $M$ is a subgroup of the complex numbers $\mathbb{C}$, there
is a pairing $< , >_\Gamma$ on $\cuco(\Tp,M)$,
given by
$$
<\varphi_1,\varphi_2>_\Gamma:=
\sum_{e\in\vec{\EE}(\Gamma\backslash\Tp)}
\varphi_1(e)\overline{\varphi_2(e)}\frac{1}{\#\overline{\overline{\Gamma}}_e},
$$
where 
$$
\overline{\overline{\Gamma}}_e := {\rm Stab}_{\overline{\overline{\Gamma}}}(e)
$$ 
is the (finite) stabiliser of $e$ in $\overline{\overline{\Gamma}}$. 

\begin{Def} \label{trace}
Let $\Gamma'$ be a subgroup of $\Gamma$ of finite index and
$\{\gamma_1,\dots,\gamma_r\}$ be a set of representatives for
$\Gamma'\backslash\Gamma$. We define the following trace map:
$$
\begin{array}{rccl}
{\rm Tr}^{\Gamma'}_\Gamma: & \cuco(\Tp,M)^{\Gamma'}
& \longrightarrow & \cuco(\Tp,M)^\Gamma \\
& \varphi & \longmapsto & \sum_{i=1}^r \gamma_i^{-1}\ast\varphi.
\end{array}
$$
\end{Def}

\begin{Lemma} \label{petersson}
The pairing above is compatible with the trace map, i.e.:
$$
<\varphi,\psi>_{\Gamma'}=<{\rm Tr}_\Gamma^{\Gamma'}\varphi,\psi>_\Gamma,
$$
where $\varphi\in\cuco(\Tp,M)^{\Gamma'}$ and
$\psi\in\cuco(\Tp,M)^\Gamma$.
\end{Lemma}

\begin{proof}
By definition,
\begin{eqnarray*}
<{\rm Tr}_\Gamma^{\Gamma'}\varphi,\psi>_\Gamma
& = & \sum_{i=1}^r <\gamma_i^{-1}\ast\varphi,\psi> \\
& = & \sum_{i=1}^r \sum_{e\in\vec{\EE}(\Gamma\backslash\Tp)}
      \varphi(\gamma_i e)\overline{\psi(e)}
      \frac{1}{\#\overline{\overline{\Gamma}}_e}.
\end{eqnarray*}
The $\Gamma$-orbit of $e$ decomposes into $\Gamma'$-orbits as follows:
$$
\Gamma e = \Gamma'\gamma_1 e \cup\ldots\cup \Gamma'\gamma_r e.
$$
This decomposition is not necessarily disjoint. We observe that
for all $i$, the number of $j$'s such that
$$
\Gamma'\gamma_i e = \Gamma'\gamma_j e,
$$
is equal to
$$
\frac{\#\overline{\overline{\Gamma}}_e}
     {\#\overline{\overline{\Gamma'}}_{\gamma_i e}}.
$$
This concludes the proof.
\end{proof}

\subsection{Automorphic forms}

Let $\KK\subseteq G(\OO)$ be an open subgroup
and $L$ a subfield of the complex numbers.
Later we will be interested in groups of the form
$$
\KK_0(\nn):=\left\{\left(
\begin{array}{cc} \ud{a} & \ud{b} \\ \ud{c} & \ud{d} \end{array}
\right)\in G(\OO) : \underline{c}\equiv 0 \textrm{ mod } \nn \right\}
$$
for an effective divisor $\nn$ of $F$.

\begin{Def}
The space of automorphic cusp forms $W_\mathfrak{p}(\KK,L)$ at $\mathfrak{p}$ consists
of functions
$$
f:G(\Ac)\longrightarrow L
$$
such that
\begin{enumerate}
\item[(i)]
for all $\gamma\in G(F)$, $\ud{g}\in G(\Ac)$ and $\ud{k}\in\KK Z(F_\mathfrak{p})$,
$$
f(\gamma\ud{g}\ud{k})=f(\ud{g})
$$
and
\item[(ii)]
for all $\ud{g}\in G(\Ac)$,
$$
\int_{F\backslash\Ac}f(\left(
\begin{array}{cc} 1 & \ud{x} \\ 0 & 1 \end{array} \right)
\ud{g}) d\ud{x}=0,
$$
where $d\underline{x}$ is a Haar measure on $F\backslash\Ac$.
\end{enumerate}
\end{Def}

\noindent
In particular, an automorphic cusp form is a function on
$$
Y_\fp(\KK) := G(F)\backslash G(\Ac)/\KK\cdot Z(F_\mathfrak{p}).
$$
There is a natural $G(F_\mathfrak{p})$-action (from the right) on the space
$$
V_{\rsp,\mathfrak{p}}(L):=\{f:\PP^1(F_\mathfrak{p})\rightarrow L:
f \textrm{ locally constant}\}/L.
$$
This is called the \textit{(L-valued) special representation} 
$\varrho_{\rsp,\mathfrak{p}}$ of $G(F_\mathfrak{p})$.
An automorphic cusp form $f\in W_\mathfrak{p}(\KK,L)$
\textit{transforms like} $\varrho_{\rsp,\mathfrak{p}}$ if the $L$-vector space 
generated
by its right $G(F_\mathfrak{p})$-translates is isomorphic to a finite number of
copies of $\varrho_{\rsp,\mathfrak{p}}$. We denote the space of such forms by
$W_{\rsp,\mathfrak{p}}(\KK,L)$.

\medskip
Assume that $\KK$ decomposes as a product
$$
\KK=\KK_{f,\mathfrak{p}}\times\II_\mathfrak{p},
$$
where $\II_\mathfrak{p}$ is the Iwahori group. Then we can choose a system
$R_\fp$ of representatives for the finite set
$G(F)\backslash G(\Ac_{f,\mathfrak{p}})/\KK_{f,\mathfrak{p}}$ and we define 
$$
\Gamma_{\underline{x}}:=G(F)\cap\underline{x}\KK_{f,\mathfrak{p}}\underline{x}^{-1}
$$
to be the intersection in $G(\Ac_{f,\mathfrak{p}})$ for $\underline{x}\in R_\fp$.
Every element $\underline{g}\in G(\Ac)$ can be written as a product
$$
\underline{g} = \gamma(\underline{x}\times 1_\fp)(\underline{k}\times 1_\fp)
(\underline{1}_{f,\fp}\times g_\fp),
$$
for some $\gamma\in G(F)\subseteq G(\Ac)$, $\underline{k}\in\KK_{f,\fp}$,
$g_\fp\in G(F_\fp)$, $1_\fp\in G(F_\fp)$ and $\underline{1}_{f,\fp}\in G(\Ac_{f,\fp})$
the respective units, and a uniquely determined $\underline{x}\in R_\fp$.
This leads to the following identification (\cite[4.5.4]{gr}):
$$
\begin{array}{rccc}
\Phi_\fp: & Y_\fp(\KK) & \stackrel{\cong}{\longrightarrow} & \displaystyle
\coprod_{\underline{x}\in R_\fp}
\Gamma_{\underline{x}}\backslash G(F_\fp)/Z(F_\fp)\cdot\II_\fp \\
& [ \underline{g} ]  & \longmapsto & [g_\fp]
,
\end{array}
$$
where $[\cdot]$ denotes the double class of an element. Of course,
the group on the right is the same as
$\coprod_{\underline{x}\in R_\fp}
\vec{\EE}(\Gamma_{\underline{x}}\backslash\Tp)$
(see sect.\ 1.1).
The following important theorem can e.g.\ be found in \cite[4.7.6]{gr}:

\begin{Thm}[Drinfeld] \label{drinfeld}
Assume $\KK=\KK_{f,\mathfrak{p}}\times\II_\mathfrak{p}$ with an open subgroup
$\KK_{f,\mathfrak{p}}\subseteq G(\OO_{f,\mathfrak{p}})$. Under the identification above,
the following spaces are isomorphic:
$$
W_{\rsp,\mathfrak{p}}(\KK,L)
\cong
\bigoplus_{\underline{x}\in R_\fp}\cuco(\Tp,L)^{\Gamma_{\underline{x}}}.
$$
\end{Thm}

\medskip\noindent
One defines a theory of Hecke operators on cuspidal harmonic cochains
which is compatible with the Hecke algebra on automorphic forms. 
Using the Petersson product, it
is therefore possible to define new- and oldforms (see below).

\begin{Rem}
Let $F=\FQ(T)$, $\deg(\mathfrak{p})=1$ and 
$\KK=\KK_0(\fp\nn)=\KK_{f,\fp}\times\II_\mathfrak{p}$, where
$\nn$ is a positive divisor such that $\mathfrak{p}{\hbox{$\not|\,$}}\nn$.
Then $R_\fp$ consists only of one element and
$$
W_{\rsp,\mathfrak{p}}(\KK_0(\mathfrak{p}\nn),L)\cong\cuco(\Tp,L)^{\Gamma_0^\fp(\nn)}.
$$
\end{Rem}

The space of automorphic cusp forms admits a non-degenerate pairing,
the so-called Petersson product, given by
$$
(f_1,f_2)_\mu:=\int_{Y_\fp(\KK)}f_1(\underline{g})\overline{f_2(\underline{g})}
d\mu(\underline{g}),
$$
where $\mu$ is a suitably normalised Haar measure on 
$G(\Ac)/Z(G(F_\infty))$.

With the identifications of theorem \ref{drinfeld}, the Haar measure
becomes
$$
\mu(e)=\frac{1}{\#(\overline{\overline{\Gamma}}_{\underline{x}})_e},
$$
and the pairing in section 1.1 is seen to be the Petersson product.
(More details can be found in \cite[4.8]{gr}.)

\subsection{Theta functions and Tate curves}

Let $E/F$ be an elliptic curve of conductor $\mathfrak{p}\nn$
with split multiplicative reduction at $\mathfrak{p}$ and $\Gamma$
an arithmetic group.
Then $E$ is a Tate curve at $\mathfrak{p}$, i.e.
$$
E(F_\mathfrak{p})\cong F_\mathfrak{p}^\times/q_{E,\mathfrak{p}}^\ZZ,
$$
where $q_{E,\fp}\in F_\fp^\times$ is the Tate period. This isomorphism can be
constructed explicitly by means of theta functions. The contents of this
section are a review of \cite[Sections 5,7 and 9]{gr}.

\begin{Thm}[Gekeler-Reversat]
Let $\omega\in\Omega_\mathfrak{p}$ be a randomly chosen point. For $\alpha\in\Gamma$,
one defines
$$
u_\alpha:\left\{
\begin{array}{ccc}
\Omega^*_\mathfrak{p} & \longrightarrow & \CC_\mathfrak{p} \\
z & \longmapsto & \displaystyle\prod_{\varepsilon\in\overline{\overline{\Gamma}}}
                     \frac{z-\varepsilon\omega}{z-\varepsilon\alpha\omega}
\end{array}
\right..
$$
This product converges to a $\CC_\fp^\times$-valued holomorphic theta function,
i.e.\ for all $\beta\in\Gamma$, $u_\alpha$
satisfies a functional equation
$$
u_\alpha(\beta z)=c_\alpha(\beta)u_\alpha(z).
$$
The definition of $u_\alpha$ is independent of the chosen
$\omega$ and only depends on the class of $\alpha$ in $\overline{\Gamma}$.
This induces a group homomorphism
$$
\overline{c}:\left\{
\begin{array}{ccc}
\overline{\Gamma} & \longrightarrow & 
                {\rm Hom}(\overline{\Gamma},F_\mathfrak{p}^\times) \\
\alpha & \longmapsto & (\beta\mapsto c_\alpha(\beta))
\end{array}
\right..
$$
Furthermore, the map
$$
\begin{array}{ccc}
\overline{\Gamma}\times\overline{\Gamma} & \longrightarrow & F_\mathfrak{p}^\times \\
(\alpha,\beta) & \longmapsto & c_\alpha(\beta)
\end{array}
$$
defines a symmetric bilinear pairing.
\end{Thm}
\noindent
We observe that by definition, $u_\gamma(\infty)=1$ for all $\gamma\in\Gamma$.

\bigskip
The elliptic curve $E$ has a uniformisation
$$
\overline{M}_0(\nn)\longrightarrow E
$$
by a Drinfeld modular curve, defined over $F$. The affine algebraic
curve $M_0(\nn)$ is the moduli scheme representing the functor of
Drinfeld modules of rank $2$ with level $\nn$-structure.
The rigid analytic variety $M_0^{\rm an}(\nn)$ associated to
$M_0(\nn)\otimes_F \CC_\fp$ decomposes as:
$$
M_0^{\rm an}(\nn)=
\coprod_{\underline{x}\in R_\fp}\Gamma_{\underline{x}}\backslash\Omega_\fp.
$$
This curve can be compactified to a projective curve $\overline{M}_0(\nn)$
by adding a finite number of points (the so-called cusps), which can be 
seen as coming from $\PP^1(F)\subset \partial\Omega_\fp$.
As above, we get a decomposition:
$$
\overline{M}_0^{\rm an}(\nn)=
\coprod_{\underline{x}\in R_\fp}\overline{M}_{\Gamma_{\underline{x}}}(\CC_\fp)=
\coprod_{\underline{x}\in R_\fp}\Gamma_{\underline{x}}\backslash\Omega_{\fp}^*,
$$
and its Jacobian decomposes accordingly.
It is enough to work on a chosen component (\cite[(9.6)]{gr}).

Assume that $E$ is the strong Weil curve in its isogeny class.
It corresponds to a primitive Hecke eigenform 
$$
c_{\mathfrak{p}}=(c_{\fp,\underline{x}})_{\underline{x}\in R_\fp}\in
\bigoplus_{\underline{x}\in R_\fp}
\cuco^{{\rm new}}(\Tp,\ZZ)^{\Gamma_{\underline{x}}}
$$ 
with rational eigenvalues. I.e.\ $c_{\mathfrak{p}}$ is normalised such
that for all $\underline{x}\in R_\fp$, 
$c_{\mathfrak{p},\underline{x}}\in j(\overline{\Gamma}_{\underline{x}})$ but
$c_{\mathfrak{p},\underline{x}}\not\in nj(\overline{\Gamma}_{\underline{x}})$ for $n>1$
(\cite[(9.1)]{gr}). 
We now choose a component $\Gamma_{\underline{x}}\backslash\Omega_\fp^\ast$,
as above, and consider the corresponding newform $c_\fp=c_{\fp,\underline{x}}$.
Let $\gamma_\mathfrak{p}\in\overline{\Gamma}_{\underline{x}}$
be its preimage under $j$. 
Then the analytic
uniformisation of $E$ at $\mathfrak{p}$ is given by the diagram:
$$
\xymatrix{
1 \ar[r] & \overline{\Gamma}_{\underline{x}} 
           \ar[r]^{\overline{c}\textrm{\ \ \ \ \ \ \ \ \ \ }} \ar[d] & 
{\rm Hom}(\overline{\Gamma}_{\underline{x}},\CC_\fp^\times) \ar[r] \ar[d]^{{\rm ev}_\fp} &
{\rm Jac}(\overline{M}_{\Gamma_{\underline{x}}}(\CC_\mathfrak{p}))
                      \ar[r] \ar[d]^{{\rm pr}_\fp} & 0 \\
1 \ar[r] & \Lambda_\mathfrak{p} \ar[r] & \CC_\mathfrak{p}^\times \ar[r] & 
E(\CC_\mathfrak{p}) \ar[r] & 0,
}
$$
where ${\rm ev}_\fp$ is the evaluation map at $\gamma_\fp$, $\Lambda_\mathfrak{p}$
is the image of $\overline{\Gamma}_{\underline{x}}$ under ${\rm ev}_\fp$ and
${\rm pr}_\fp$ is the map induced by ${\rm ev}_\fp$.
It is explicitly given by:
$$
\begin{array}{rccc}
{\rm pr}_\fp: & {\rm Jac}(\overline{M}_{\Gamma_{\underline{x}}}(\CC_\fp)) &
    \longrightarrow & E(\CC_\fp) \\
& (a)-(b) & \longmapsto & \frac{u_{\gamma_\fp}(a)}{u_{\gamma_\fp}(b)}
\end{array}
$$
where $a,b\in\Omega_\fp^\ast$.

\pagebreak
\noindent
There exists a divisor $d$ of $q^{\deg(\mathfrak{p})}-1$ and $t\in F_\mathfrak{p}^\times$
with $|t|_\fp<1$ (both dependent on $c_\mathfrak{p}$) such that
$\Lambda_\mathfrak{p}=\boldsymbol\mu_d\times t^\ZZ$ 
($\boldsymbol\mu_d$ the $d$-th roots of unity in $F_\mathfrak{p}^\times$).
The Tate period is given by $q_{E,\mathfrak{p}}=t^d$.

In the case that $A_\mathfrak{p}$ is isomorphic to a polynomial ring,
Gekeler shows in \cite{gek} that $d=1$. In particular,
$$
\Lambda_\mathfrak{p}=q_{E,\mathfrak{p}}^\ZZ
$$
if $F=\FQ(T)$ and $\mathfrak{p}$ is of degree $1$.

\bigskip
An analogue of the well known Manin-Drinfeld theorem was established
in \cite{gek3}:

\begin{Thm}[Gekeler] \label{mandri}
For each congruence subgroup $\Gamma_0^\fp(\nn)$ of $\gl_2(F)$, 
the cuspidal divisor group
\begin{eqnarray*}
\Cc_\Gamma & := &
{\rm Div}_0(\overline{M}_0(\nn)(\CC_\fp)-M_0(\nn)(\CC_\fp))
\mod ({\rm principal\ divisors}) \\
& = & {\rm Div}_0(\Gamma_0^\fp(\nn)\backslash\PP^1(F))
\mod ({\rm principal\ divisors})
\end{eqnarray*}
is finite.
\end{Thm}

\subsection{Multiplicative integrals}

Based on one of Bertolini's ideas,
the following multiplicative integrals are considered in \cite{l}.
Let $X$ be a topological space such that a basis for its topology is given
by its open compact subsets.
Given a continuous function $f:X\rightarrow \CC^\times_\fp$ and a
measure $\mu\in{\rm Meas}(X,\ZZ)$, i.e.\ a finitely additive function
on the compact open subsets of $X$, we define:
$$
\mint_X f(t)d\mu(t):=\varinjlim_\alpha\prod_{U\in\Cc_\alpha}
f(t)^{\mu(U)},
$$
where $\{\Cc_\alpha\}$ is the direct system of finite covers of $X$
by compact open subsets, and $t\in U$ is chosen arbitrarily. This
integral exists and is independent of the choices of the $t$'s.
Furthermore, it induces a continuous homomorphism from the group
of continuous functions $\Cc(X,\CC_\fp^\times)$ to $\CC_\fp^\times$.

For future reference, we state one important (but obvious) property:

\begin{Lemma} \label{intprop1}
Let $X$ and $\mu$ be as above. Then
$$
\mint_X cd\mu(t) = c^{\mu(X)},
$$
for all $c\in\CC_\fp^\times$.
\end{Lemma}

If we choose $X$ to be the boundary $\partial\Omega_\fp$, the computation
of a multiplicative integral can be accomplished as follows.
Choose a vertex $v\in\vec{\EE}(\Tp)$ and, for $e\in\vec{\EE}(\Tp)$
pointing away from $v$, define $\textrm{dist}_v(e)$
to be the distance between $o(e)$ and $v$ (i.e.\ the number of
edges of a geodesic). Then, for 
$f\in\Cc(\partial\Omega_\fp,\CC_\fp^\times)$,
$$
\mint_{\partial\Omega_\mathfrak{p}} f(t)d\mu(t) = 
\lim_{n\rightarrow\infty}\prod_{{\rm dist}_v(e)=n}
f(t_e)^{\mu(U(e))},
$$
where $t_e\in U(e)$ is chosen arbitrarily. 
This is independent of the choice of $v$, which is normally
chosen to be the base vertex $v_0$ (see sect.\ 1.1).

Another important property is the following:

\begin{Lemma} \label{intprop2}
Let $U$ be an open subset of $\partial\Omega_\fp$ and
$\gamma\in\gl_2(F_\fp)$. Then
$$
\mint_{\gamma U}fd(\gamma\ast\mu) = \mint_U(f\circ\gamma)d\mu,
$$
where $\gamma\ast\mu(V)=\mu(\gamma^{-1}V)$.
\end{Lemma}

It is possible to obtain a multiplicative version of Teitelbaum's
Poisson inversion formula (\cite{t2}):

\begin{Prop} \label{poisson}
Let $\alpha\in\GaB$. Then
$$
\frac{u_\alpha(z_2)}{u_\alpha(z_1)} =
\mint_{\partial\Omega_\mathfrak{p}}\frac{t-z_2}{t-z_1}d\mu_{j(\alpha)}(t),
$$
for all $z_1,z_2\in\Omega^*_\mathfrak{p}$.
\end{Prop}

\begin{proof}
This is a combination of propositions 8 and 24 in \cite{l}.
\end{proof}

%% file: integr.tex
\section{Integration on $\Tp\times\Tinf/\Gamma$}

\subsection{Harmonic cochains}

We choose a further place $\infty$ of $F$, different of $\fp$, and set 
$$
S:=\{\fp,\infty\}.
$$
We define the following subgroup of $\textrm{GL}_2(A_S)$:
$$
\Gamma:=\left\{\gamma=
\left(\begin{array}{cc}a & b \\ c & d \end{array}\right)
\in\textrm{GL}_2(A_S):
\det(\gamma)\in\FQ \text{\rm\ and } c\in\nn
\right\}
$$
where $\nn\leq A_S$ is an ideal corresponding to an effective divisor of $F$, 
which we will also call $\nn$, s.t.\ $\mathfrak{p},\infty{\hbox{$\not|\,$}}\nn$.
We will use the same notation for the corresponding ideals in $A_\mathfrak{p}$ and
$A_\infty$.

\begin{Lemma} \label{density}
Let $\fq\in S$. Then
the group $\Gamma\cap\SL_2(F_\fq)$ is dense in $\SL_2(F_\fq)$.
\end{Lemma}

\begin{proof}
Let $\Gamma^c$ be the closure of $\Gamma\cap\SL_2(F_\fq)$ in $\SL_2(F_\fq)$.
By the strong approximation theorem, $A_S$ and $\nn$ are dense in $F_\fq$. 
Therefore, since $\Gamma$ contains the subgroups 
$$
\begin{pmatrix} 1 & 0 \\ \nn & 1 \end{pmatrix}
\textrm{\ \ \ \ and\ \ \ \ }
\begin{pmatrix} 1 & A_S \\ 0 & 1 \end{pmatrix},
$$
we see that $\Gamma^c$ contains the full triangular subgroups 
$$
\begin{pmatrix} 1 & 0 \\ F_\fq & 1 \end{pmatrix}
\textrm{\ \ \ \ and\ \ \ \ }
\begin{pmatrix} 1 & F_\fq \\ 0 & 1 \end{pmatrix}.
$$
Let $\alpha\in\SL_2(F_\fq)$ and put $a:=\alpha(\infty)$, $b:=\alpha(0)$. 
The matrix 
$$
\beta := \begin{pmatrix} 1 & b \\ 0 & 1 \end{pmatrix}
\begin{pmatrix} 1 & 0 \\ (a-b)^{-1} & 1 \end{pmatrix}\in\Gamma^c
$$
also sends $\infty$ to $a$ and $0$ to $b$. It follows that $\beta^{-1}\alpha$ is a diagonal 
matrix in $SL_2(F_\fq)$, i.e., of the form
$$
\begin{pmatrix} x & 0 \\ 0 & x^{-1} \end{pmatrix}.
$$
But all such diagonal matrices are contained in $\Gamma^c$, since
$$
\begin{pmatrix} x & 0 \\ 0 & x^{-1} \end{pmatrix}
= \begin{pmatrix} 1 & 0 \\ x^{-1}-1 & 1 \end{pmatrix}
\begin{pmatrix} 1 & 1 \\ 0 & 1 \end{pmatrix}
\begin{pmatrix} 1 & 0 \\ x-1 & 1 \end{pmatrix}
\begin{pmatrix} 1 & -x^{-1} \\ 0 & 1 \end{pmatrix}
$$
(this is a special case of Whitehead's lemma).
\end{proof}

In the following we use the letters $w$ and $s$ to denote vertices
and edges on $\Tinf$, whereas we keep $v$ and $e$ for $\Tp$.

\begin{Cor}
The group $\Gamma$ acts transitively on the sets of unoriented edges of
$\Tp$ and $\Tinf$. Its action on oriented edges has two orbits:
$$
\vec{\EE}(\Tp)=\Gamma e_0\coprod\Gamma\overline{e}_0
\textrm{\ \ \ \ and\ \ \ \ }
\vec{\EE}(\Tinf)=\Gamma s_0\coprod\Gamma\overline{s}_0.
$$
\end{Cor}

\begin{proof}
We observe that $\vec{\EE}(\Tp)=\SL_2(F_\fp) e_0\coprod\SL_2(F_\fp)\overline{e}_0$.
Let $e$ be in the $\SL_2(F_\fp)$-orbit of $e_0$. The set
$$
T_e:=\{\gamma\in\SL_2(F_\fp):\gamma e_0=e\}
$$
is open in $\SL_2(F_\fp)$, since it is a coset of the stabiliser of $e_0$.
Lemma \ref{density} implies that
$$
T_e\cap\Gamma\not=\emptyset,
$$
i.e., $e$ is in the $\Gamma$-orbit of $e_0$. We can apply the same argument
for the orbit of $\overline{e}_0$. 
The statement for $\Tinf$ follows analogously.
\end{proof}

\medskip
We want to study harmonic cochains on the product
$\Tree=\Tp\times\Tinf$. Let $\vec{\EE}(\Tree)$ (resp.\ $\VV(\Tree)$) denote
the product of the corresponding sets of (directed) edges (resp.\ vertices).
We observe that for the base vertex $v_0\in\VV(\Tp)$,
$$
\begin{array}{rcl}
\Gamma_{v_0}:=\Gamma\cap\mbox{Stab}(v_0)
& = & \Gamma_0^\infty(\nn) \\
\end{array}
$$
and for the base edge $e_0\in\VV(\Tp)$,
$$
\begin{array}{rcl}
\Gamma_{e_0}:=\Gamma\cap\mbox{Stab}(e_0)
& = & \Gamma_0^\infty(\mathfrak{p}\nn).
\end{array}
$$
The stabilisers of the other vertices in $\Tp$ are $\Gamma$-conjugate 
to $\Gamma_{v_0}$. For an edge $e\in\vec{\EE}(\Tp)$, either
the stabiliser of $e$ or $\overline{e}$ is conjugate to 
$\Gamma_{e_0}$, depending on whether $e$ is in the $\Gamma$-orbit
of $e_0$ or $\overline{e}_0$.
The same holds for the base vertex/edge $w_0,s_0$
on $\Tinf$ with the corresponding subgroups of $\mbox{GL}_2({A_\mathfrak{p}})$.

\begin{Def} \label{codef}
Let $M$ be an abelian group.
The space of $\Gamma$-invariant harmonic cochains $\cuco(\Tree,M)^\Gamma$
with values in $M$, cuspidal with respect to $S$, 
consists of functions
$c:\vec{\EE}(\Tree)\rightarrow M$ such that:
\begin{enumerate}
\item[(i)] $\forall\gamma\in\Gamma\textrm{ }
\forall (e,s)\in\vec{\EE}(\Tree): 
c(\gamma e,\gamma s)=c(e,s)$\\
and $c$ has compact support modulo $\Gamma$,
\item[(ii)] $\forall e\in\vec{\EE}(\Tp):c(e,\cdot)\in\co(\Tinf,M)$ and
\item[(iii)] $\forall s\in\vec{\EE}(\Tinf):
c(\cdot,s)\in\co(\Tp,M)$.
\end{enumerate}
\end{Def}

\bigskip
\noindent\textbf{Remark.}
These conditions imply that
\begin{enumerate}
\item[]
$\forall e\in\vec{\EE}(\Tp):c(e,\cdot)\in\hc{\infty}{M}{\Gamma_e}$
and
\item[]
$\forall s\in\vec{\EE}(\Tinf): c(\cdot,s)\in\hc{\mathfrak{p}}{M}{\Gamma_s}$.
\end{enumerate}

\begin{Def}
Let $L$ be a subgroup of the complex numbers.
The subspace of $\mathfrak{p}$-newforms
$\cuco^{\mathfrak{p}\text{\rm -new}}(\Tinf,L)^{\Gamma_0^\infty(\mathfrak{p}\nn)}$ 
of $\cuco(\Tinf,L)^{\Gamma_0^\infty(\fp\nn)}$ is defined
to be the orthogonal complement of the images of the natural inclusions of
$\cuco(\Tinf,M)^{\Gamma_{v_0}}$ and $\cuco(\Tinf,M)^{\Gamma_{v_{-1}}}$ 
into $\cuco(\Tinf,M)^{\Gamma_{e_0}}$, induced by the equality
$$
\Gamma_{e_0}=\Gamma_{v_0}\cap\Gamma_{v_{-1}},
$$
with respect to the Petersson product.
\end{Def}

The following is an analogue of \cite[Lemma 1.3,3.]{da}:

\begin{Prop} \label{newiso}
The map
$$
\begin{array}{ccc}
\cuco(\Tree,L)^\Gamma & \longrightarrow &
      \cuco^{\mathfrak{p}\text{\rm -new}}(\Tinf,L)^{\Gamma_0^\infty(\mathfrak{p}\nn)} \\
c & \longmapsto & c_\infty(\cdot):=c(e_0,\cdot)
\end{array}
$$
is a well-defined isomorphism.
\end{Prop}

\begin{proof} 
Since $\vec{\EE}(\Tp)=\Gamma e_0\coprod\Gamma\overline e_0$, $c$ is completely 
determined by $c(e_0,\cdot)$. On the other hand, given 
$c_\infty\in\hc{\infty}{L}{\Gamma_0(\mathfrak{p}\nn)}$, we can define a cochain 
$c\in\cuco(\Tree,L)^\Gamma$ by
$c(\gamma e_0,s):=c_\infty(\gamma^{-1}s)$ and
$c(\gamma \overline{e}_0,s):=-c_\infty(\gamma^{-1}s)$ for $\gamma\in\Gamma$,
respectively, provided the condition of harmonicity at the
vertices $v_0$ and $v_{-1}$ is satisfied (observe that there are also only two
$\Gamma$-orbits for vertices on $\Tp$). 

\noindent
Let $\{\gamma_1,\dots,\gamma_r\}$ ($r=q^{\deg(\mathfrak{p})}+1$) 
be a set of representatives of the quotient
$\Gamma_0^\infty(\nn)/\Gamma_0^\infty(\mathfrak{p}\nn)$.
Then the map
$$
\gamma_i\longmapsto\gamma_i v_{-1}
$$
defines a bijection between 
$\Gamma_0^\infty(\nn)/\Gamma_0^\infty(\mathfrak{p}\nn)$
and the $\Gamma_{v_0}$-orbit of $v_{-1}$, i.e.\ the set of edges with origin $v_0$.
Condition (iii) above yields 
$$
0=\sum_{o(e)=v_0}c(e,s)=
\sum_{i=1}^r c(\gamma_i\overline{e}_0,s)=
-\sum_{i=1}^r c_\infty(\gamma_i^{-1}s).
$$ 
I.e., $c_\infty$ is in the kernel of the trace map
$$
\begin{array}{rccc}
{\rm Tr}_{\Gamma_0^\infty(\nn)}^{\Gamma_0^\infty(\mathfrak{p}\nn)}:
     & \hc{\infty}{L}{\Gamma_0^\infty(\mathfrak{p}\nn)} & \longrightarrow & 
       \hc{\infty}{L}{\Gamma_0^\infty(\nn)} \\
& \psi & \longmapsto & \sum\psi\circ\gamma_i^{-1}=\sum\gamma_i\ast\psi
\end{array}
$$
By lemma \ref{petersson},
$$
<\varphi,\psi>_{\Gamma_0^\infty(\mathfrak{p}\nn)}=
<{\rm Tr}_{\Gamma_0^\infty(\nn)}^{\Gamma_0^\infty(\mathfrak{p}\nn)}\varphi,\psi>_{\Gamma_0^\infty(\nn)}
$$
for $\varphi\in\hc{\infty}{L}{\Gamma_0^\infty(\mathfrak{p}\nn)}$ and 
$\psi\in\hc{\infty}{L}{\Gamma_0^\infty(\nn)}$.

\noindent
The edges exiting from $v_{-1}$ are the orbit of $e_0$ in $\Gamma_{v_{-1}}$.
Similar reasoning as above implies that 
$$
Tr^{\Gamma_0^\infty(\mathfrak{p}\nn)}_{\Gamma_{v_{-1}}}(c_\infty)=0.
$$ 
Hence $c_\infty$ is orthogonal to 
$\hc{\infty}{L}{\Gamma_{v_{-1}}}$. 
\end{proof}

Now we consider the following group:
$$
\tilde{\Gamma}:=
\left\{\gamma=
\left(\begin{array}{cc}a & b \\ c & d \end{array}\right)
\in\textrm{GL}_2(A_S): c\in\nn
\right\}.
$$
We observe that
$$
\Gamtil_{v_0}\cap\gl_2(A_\infty)=\Gamma_0^\infty(\nn)
$$
and
$$
\Gamtil_{e_0}\cap\gl_2(A_\infty)=\Gamma_0^\infty(\fp\nn),
$$
where $\Gamtil_\cdot:=\textrm{Stab}_{\Gamtil}(\cdot)$.

\medskip
We define the space of
$\Gamtil$-invariant, $S$-cuspidal harmonic cochains
$$\cuco(\Tree,M)^\Gamtil$$ with values in an abelian group $M$
as in definition \ref{codef}.
This is a subspace of $\cuco(\Tree,M)^\Gamma$.

\medskip
By the $S$-unit theorem, there exists $u\in F$ such that
$$
A_S^\times=\FQ^\times\times u^\ZZ
$$
and
$$
\textrm{div}(u)=h_\fp\cdot[\fp]-h_\infty\cdot[\infty].
$$
The group $\Gamtil$ is generated by $\Gamma$ and the matrix
$$
\gamma=
\begin{pmatrix} u & 0 \\ 0 & 1 \end{pmatrix}
\in \Gamtil-\Gamma.
$$
In the case that $h_\fp$ is odd, $\gamma e_0$ is in the $\Gamma$-orbit
of $\overline{e}_0$ and there exists $\beta\in\Gamma$ such that
$\gamma e_0=\beta\overline{e}_0$. Alternatively, if $h_\fp$ is even,
there exists $\beta\in\Gamma$ such that $\gamma e_0=\beta e_0$.
In either case we define
$$
\alpha_\fp:=\beta^{-1}\gamma\in\Gamtil-\Gamma.
$$ 
Then the Atkin-Lehner involution at $\fp$
on $\cuco(\Tinf,L)^{\Gamma_0^\infty(\fp\nn)}$ is given by:
$$
W_\fp c_\infty:=\alpha_\fp\ast c_\infty.
$$
An easy calculation shows that the image of $\cuco(\Tree,L)^\Gamtil$ in
$\cuco^{\fp\textrm{-new}}(\Tinf,L)^{\Gamma_0^\infty(\fp\nn)}$
(via the isomorphism in proposition \ref{newiso}) is
$$
\cuco^{\fp{\rm-new},-}(\Tinf,L)^{\Gamma_0^\infty(\fp\nn)},
$$
the space of those $\fp$-newforms which are eigenforms under $W_\fp$
of eigenvalue $(-1)^{h_\fp}$.
This shows:

\begin{Cor} \label{cornewiso}
The following diagram (with horizontal isomorphisms) is commutative:
$$
\begin{array}{ccc}
\cuco(\Tree,M)^\Gamtil & \stackrel{\cong}{\longrightarrow} &
\cuco^{\fp{\rm-new},-}(\Tinf,L)^{\Gamma_0^\infty(\fp\nn)} \\
\bigcap && \bigcap \\
\cuco(\Tree,M)^\Gamma & \stackrel{\cong}{\longrightarrow} &
\cuco^{\fp{\rm -new}}(\Tinf,L)^{\Gamma_0^\infty(\fp\nn)}.
\end{array}
$$
\end{Cor}

\bigskip
\noindent\textbf{Remark.}
Similarly, by interchanging $\mathfrak{p}$ and $\infty$ 
and exchanging $e_0$ with $s_0$ we get:
$$
\begin{array}{ccc}
\cuco(\Tree,L)^\Gamtil & \stackrel{\cong}{\longrightarrow} &
\cuco^{\infty{\rm -new},-}(\Tp,L)^{\Gamma_0^\mathfrak{p}(\nn\infty)} \\
\bigcap && \bigcap \\
\cuco(\Tree,L)^\Gamma & \stackrel{\cong}{\longrightarrow} &
\cuco^{\infty{\rm -new}}(\Tp,L)^{\Gamma_0^\mathfrak{p}(\nn\infty)} \\
c & \longmapsto & c_\mathfrak{p}(\cdot).
\end{array}
$$

\subsection{Double integrals}

Let $x,y\in\PP^1(F)$. We denote the geodesic joining the corresponding
ends on $\Tinf$ by $A(x,y)$.
For any $c\in\cuco(\Tree,M)^\Gamtil$, it follows from the definitions that
$$
c\{x\rightarrow y\}:
e\longmapsto\sum_{s\in A(x,y)}c(e,s)
$$
is a $M$-valued harmonic cochain on $\Tp$. (Note that since
$c$ is cuspidal, this is a finite sum for all $e$.)
Therefore it corresponds to a measure
$$
\mu_c\{x\rightarrow y\}
$$
on $\PP^1(F_\mathfrak{p})$ of total mass $0$.

\medskip
If we choose $\gamma\in\Gamtil$ such that $e=\gamma e_0$, we see that
\begin{eqnarray*}
c\{x\rightarrow y\}(e)& = & 
\sum_{s\in A(\gamma^{-1}x,\gamma^{-1}y)} c_\infty(s) \\
& = & \sum_{s\in A(x,y)} \gamma\ast c_\infty(s)
\end{eqnarray*}
(compare cor.\ \ref{cornewiso}). Or equivalently,
$$
\mu_c\{x\rightarrow y\}(U(e))
= -\sum_{s\in A(x,y)} \gamma\ast c_\infty(s).
$$
This implies the following:
\begin{Lemma} \label{g-inv} 
For all $c\in\cuco(\Tree,M)^\Gamtil$, the map 
$(x,y)\mapsto\mu_c\{x\rightarrow y\}$ 
induces a $\Gamtil$-module homomorphism 
$$
\mathcal{M}\longrightarrow {\rm Meas}_0(\partial\Omega_\mathfrak{p},M)
$$
from the group of modular symbols $\mathcal{M}={\rm Div}^0(\PP^1(F))$
to measures on $\partial\Omega_\mathfrak{p}$.
\end{Lemma}
\noindent
\textbf{Remark.} The analogous statement for $\Gamma$-invariant cochains is
also true.

\bigskip
For $c\in\cuco(\Tree,\ZZ)^\Gamtil$,
we define the following double (multiplicative) integral:
$$
\mint_{z_1}^{z_2}\int_x^y\omega :=
\mint_{\PP^1(F_\mathfrak{p})}\left( \frac{t-z_2}{t-z_1} \right) d\mu_c\{x\rightarrow y\}(t)
\in\CC^\times_\mathfrak{p}
$$
for $z_1,z_2\in\Omega_\mathfrak{p}$.

\bigskip
\textbf{Remark:}
Firstly, we observe that this definition depends on the choice of $c$.
Since later on we will work with an especially chosen $c_E$, we dropped
it from the notation.
Secondly, although $\omega$ itself is not defined, this definition
should be regarded as a period for a "rigid analytic modular form
of weight $(2,2)$ on 
$(\Omega_\mathfrak{p}\times\Omega_\infty)/\Gamma$" 
(compare \cite{da}).

\begin{Lemma} \label{int_prop}
The double integrals above have the following properties:
\begin{eqnarray*}
\mint_{z_1}^{z_3}\int_{x_1}^{x_2}\omega &=&
\mint_{z_1}^{z_2}\int_{x_1}^{x_2}\omega \times
\mint_{z_2}^{z_3}\int_{x_1}^{x_2}\omega, \\
\mint_{z_1}^{z_2}\int_{x_1}^{x_3}\omega &=&
\mint_{z_1}^{z_2}\int_{x_1}^{x_2}\omega \times
\mint_{z_1}^{z_2}\int_{x_2}^{x_3}\omega, \\
\mint_{\gamma z_1}^{\gamma z_2}\int_{\gamma x_1}^{\gamma x_2}\omega &=&
\mint_{z_1}^{z_2}\int_{x_1}^{x_2}\omega,
\end{eqnarray*}
for all $z_1,z_2,z_3\in\Omega_\mathfrak{p}$,
$x_1,x_2,x_3\in\PP^1(F)$ and $\gamma\in\tilde{\Gamma}$.
\end{Lemma}

\begin{proof}
The first two properties follow directly from the definitions, and the third
property from:
\begin{eqnarray*}
\lefteqn{
\mint_{\PP^1(F_\fp)}\frac{t-\gamma z_2}{t-\gamma z_1}
     d\mu_c\{\gamma x_1\rightarrow\gamma x_2\}(t)
} \\
& = &
\mint_{\gamma(\PP^1(F_\fp))}\frac{t-\gamma z_2}{t-\gamma z_1}
     d\gamma\ast\gamma^{-1}\ast\mu_c\{\gamma x_1\rightarrow\gamma x_2\}(t) \\
& \stackrel{\ref{intprop2}}{=} &
\mint_{\PP^1(F_\fp)}\frac{\gamma t-\gamma z_2}{\gamma t-\gamma z_1}
     d\gamma^{-1}\ast\mu_c\{\gamma x_1\rightarrow\gamma x_2\}(t) \\
& \stackrel{\ref{g-inv}}{=} &
\mint_{\PP^1(F_\fp)}
     \frac{cz_1+d}{cz_2+d}\cdot\frac{t-z_2}{t-z_1}
     d\mu_c\{x_1\rightarrow x_2\}(t) \\
& \stackrel{\ref{intprop1}}{=} &
\mint_{\PP^1(F_\fp)}
     \frac{t-z_2}{t-z_1}d\mu_c\{x_1\rightarrow x_2\}(t),
\end{eqnarray*}
where $\gamma=\begin{pmatrix} a & b \\ c & d \end{pmatrix}\in\Gamtil$. 
\end{proof}

\subsection{Automorphic cusp forms}

For an open subgroup $\KK$ of $G(\OO)$ and a field $L$ of characteristic
zero,
we denote the intersection of the spaces of automorphic cusp forms 
at $\mathfrak{p}$ and at $\infty$ in the following way:
$$
W_S(\KK,L):=
W_{\mathfrak{p}}(\KK,L)\cap W_{\infty}(\KK,L)
$$
and similarly for $W_{\rsp}$ and newforms. 

Assume that $\KK$ decomposes as
$$
\KK=\KK_{f,S}\times\II_\fp\times\II_\infty.
$$
Since $\KK_{f,S}$ is of finite index in $\OO_{f,S}$ and
$G(F)\backslash G(\Ac_{f,S})/G(\OO_{f,S})$ is finite (of cardinality
equal to the class number ${\rm cl}(A_S)$ of $A_S$, 
compare \cite[(4.1.4)]{gr}), the set
$$
G(F)\backslash G(\Ac_{f,S})/\KK_{f,S}
$$
is finite. Let $R_S$ be a system of representatives and
$$
\Gamma_{\overline{x}}:=G(F)\cap\underline{x}\KK_{f,S}\underline{x}^{-1}
\subseteq G(\Ac_{f,S}),
$$
for $\underline{x}\in R_S$. Every element $\underline{g}\in G(\Ac)$
can be written as
$$
\underline{g} = \gamma(\underline{x}\times 1_\fp\times 1_\infty)
(\underline{k}\times 1_\fp\times 1_\infty)
(\underline{1}_{f,S}\times g_\fp\times g_\infty),
$$
for some $\gamma\in G(F)$, $\underline{k}\in\KK_{f,S}$,
$g_\fp\in G(F_\fp)$, $g_\infty\in G(F_\infty)$ and a uniquely
determined $\underline{x}\in R_S$. We define
$$
Y_S(\KK) := G(F)\backslash G(\Ac)/\KK\cdot Z(F_\fp\times F_\infty)
$$
and get the following generalisation of \cite[4.5.4]{gr}.

\begin{Lemma} \label{Y_S}
The following map is a well-defined isomorphism:
$$
\begin{array}{rccc}
\Phi_S: & Y_S(\KK) & \stackrel{\cong}{\longrightarrow} & \displaystyle
\coprod_{\underline{x}\in R_S}
\Gamma_{\underline{x}}\backslash 
G(F_\fp\times F_\infty)/Z(F_\fp\times F_\infty)\cdot(\II_\fp\times\II_\infty)\\
& [ \underline{g} ]  & \longmapsto & [g_\fp\times g_\infty].
\end{array}
$$
\end{Lemma}

\begin{proof}
Straightforward.
\end{proof}

\noindent
Again, we observe that the right hand side is 
equal to
$$
\coprod_{\underline{x}\in R_S}
\vec{\EE}(\Gamma_{\underline{x}}\backslash\Tree)
$$
(as in section 1.2). 

\begin{Rem}
Since the class groups of $A_\fp$ and $A_\infty$ surject onto the class
group of $A_S$,
the class number ${\rm cl}(A_S)$ fulfils:
$$
{\rm cl}(A_S)\bigm|\gcd({\rm cl}(A_\fp),{\rm cl}(A_\infty)).
$$
In other words,
$$
\#R_S\bigm|\gcd(\#R_\fp,\#R_\infty),
$$
where $R_\fp$ and $R_\infty$ are the corresponding systems of representatives
at $\fp$ and $\infty$.
E.g., $R_S$ consists only of one element if $F$ is the rational function field
and one of the places is of degree one.
\end{Rem}

\textbf{CONVENTION:}
As already mentioned in section 1.3, the results in \cite{gr} are obtained
by restricting, without loss of generality, to a single $\underline{x}\in R_\fp$.
The same method applies in our case. To keep our notations simple, we will assume
that $R_S$ consists only of one element. To obtain the general case, one has to
slightly generalise the definitions and results in the preceding sections.

\medskip
If we choose $\KK=\KK_0(\fp\nn\infty)$, then it decomposes as above
(see remark after thm.\ \ref{drinfeld}).
Furthermore, we can easily check that
$$
G(F)\cap\KK_{f,S}=\Gamtil.
$$
The image of the natural inclusion
$
W_{{\rm sp},S}\hookrightarrow W_{{\rm sp},\fp}
$
can easily be identified as
$$
W_{{\rm sp},\fp}^{\infty{\rm-new},-}(\KK_0(\fp\nn\infty)),
$$
the space corresponding to 
$\cuco^{\infty{\rm-new},-}(\Tp,L)^{\Gamma_0^\mathfrak{p}(\nn\infty)}$ 
via the isomorphism in theorem \ref{drinfeld}.
(Again, the analogue statement is true if we interchange $\fp$ and $\infty$.)
This leads us to the following generalisation of Drinfeld's theorem:

\begin{Thm}
The following spaces are isomorphic:
$$
W_{\rsp,S}(\KK_0(\mathfrak{p}\nn\infty),L)
\cong
\cuco(\Tree,L)^\Gamtil.
$$
Furthermore, the following diagram, which consists entirely of
isomorphisms, is commutative:
$$
\xymatrix{
W_{{\rm sp},\fp}^{\infty{\rm-new},-}(\KK_0(\fp\nn\infty)) \ar[rr]^{\cong} && 
\cuco^{\infty{\rm-new},-}(\Tp,L)^{\Gamma_0^\mathfrak{p}(\nn\infty)} \\
W_{\rsp,S}(\KK_0(\fp\nn\infty),L) 
\ar[rr]^{\cong} \ar[u]^{\cong} \ar[d]_{\cong} && 
\cuco(\Tree,L)^\Gamtil \ar[u]^{\cong} \ar[d]_{\cong} \\
W_{{\rm sp},\infty}^{\fp{\rm-new},-}(\KK_0(\fp\nn\infty)) \ar[rr]^{\cong} && 
\cuco^{\fp{\rm-new},-}(\Tinf,L)^{\Gamma_0^\infty(\fp\nn)}.
}
$$
\end{Thm}

\begin{proof}
The corresponding commutative diagram of the underlying spaces is the
following:
$$
\xymatrix{
Y_\fp(\KK) \ar[rr]^{\Phi_\fp,\cong} \ar@{->>}[d] && 
\vec{\EE}(\Gamma_0^\fp(\nn\infty)\backslash\Tp) \ar@{->>}[d] & e \ar@{|->}[d]\\
Y_S(\KK) \ar[rr]^{\Phi_S,\cong} && \vec{\EE}(\Gamtil\backslash\Tree) &
{\begin{array}{c} (e,s_0) \\ (e_0,s) \end{array}} \\
Y_\infty(\KK) \ar[rr]^{\Phi_\infty,\cong} \ar@{->>}[u] &&
\vec{\EE}(\Gamma_0^\infty(\fp\nn)\backslash\Tinf) \ar@{->>}[u] & s \ar@{|->}[u]
}
$$
(the vertical maps are surjections).
The theorem follows from theorem \ref{drinfeld} and corollary
\ref{cornewiso} (resp.\ its analogue in the succeeding remark).
\end{proof}

\textbf{Remark.}
Note that this is a slightly simplified formulation of the theorem, due
to our assumption that $R_S$ contains only one element. In general, we
have to consider direct sums of spaces of harmonic cochains on the
right hand side in the diagram. E.g.,
$$
W_{\rsp,S}(\KK_0(\mathfrak{p}\nn\infty),L)
\cong
\bigoplus_{\underline{x}\in R_S}\cuco(\Tree,L)^{\Gamma_{\underline{x}}}.
$$

%% file: Da_conj.tex
\section{Teitelbaum's conjecture}

\subsection{Darmon's period}

Let $K=F\times F$ and choose an $F$-algebra embedding
$\psi :K\rightarrow \textrm{M}_2(F)$. This induces an action of
$\psi(K^\times)$ on $\Omega^*_\mathfrak{p}$ with two fixed points
$x_\psi, y_\psi\in\PP^1(F)$.
The image of $\psi(K^\times)\cap\Gamtil$ in 
${\rm PGL}_2(F)$ is of rank $1$.
We choose a generator $\gamma_\psi$ of the free part.
Assume that $x_\psi$ is repulsive and $y_\psi$ attractive w.r.t.\ 
$\gamma_\psi$.
To such an embedding and $c\in\cuco(\Tree,\ZZ)^\Gamtil$, 
we can associate a period
$$
I_\psi=I_{\psi,c}:=\mint_z^{\gamma_\psi z}\int_{x_\psi}^{y_\psi}\omega
\in\CC^\times_\mathfrak{p},
$$
where we choose $z\in\Omega_\mathfrak{p}$ arbitrarily.

\begin{Lemma}
The period $I_{\psi}$ is independent of the choice of
$z\in\Omega_\mathfrak{p}$.
\end{Lemma}

\begin{proof}
\begin{eqnarray*}
\mint_{z_1}^{\gamma_\psi z_1}\int_{x_\psi}^{y_\psi}\omega \div
\mint_{z_2}^{\gamma_\psi z_2}\int_{x_\psi}^{y_\psi}\omega & = &
\mint_{z_1}^{z_2}\int_{x_\psi}^{y_\psi}\omega \div
\mint_{\gamma_\psi z_1}^{\gamma_\psi z_2}\int_{x_\psi}^{y_\psi}\omega \\
& = &
\mint_{z_1}^{z_2}\int_{x_\psi}^{y_\psi}\omega \div
\mint_{z_1}^{z_2}\int_{x_\psi}^{y_\psi}\omega \\
& = & 1,
\end{eqnarray*}
where the second equality follows from the third property in
lemma \ref{int_prop}.
\end{proof}

\noindent\textbf{Remark.}
Note that, unlike in Darmon's setting, $I_\psi$ depends on the choice
of $\gamma_\psi$. This is due to the fact that the group of roots
of unity in $A_S^\times$ is nontrivial (see the remark after
Prop.\ \ref{I as integral}).

\bigskip
To simplify calculations, we consider only
the case where $\psi$ is the diagonal embedding. The fixed points
of $\psi$ are $\{\infty,0\}$.
As already mentioned in section 2.1, by the $S$-unit theorem
there exists $u\in F$ such that
$$
A_S^\times=\FQ^\times\times u^\ZZ
$$
and
$$
\textrm{div}(u)=h_\fp\cdot[\fp]-h_\infty\cdot[\infty].
$$
We make the following choice:
$$
\gamma_\psi=
\begin{pmatrix} u & 0 \\ 0 & 1 \end{pmatrix}
\in \Gamtil.
$$
It acts as shift by $h_\fp$ along the axis $A(\infty,0)$ on $\Tp$ and as shift
by $h_\infty$ along $A(0,\infty)$ on $\Tinf$.

\bigskip
Following Darmon, we define the winding element attached to $\psi$
as
$$
W_\psi:=\sum_{e\in A(v,\gamma_\psi v)}\mcxy(U(e)),
$$
where $v\in\VV(\Tp)$ is chosen arbitrarily. 
By observing that $\cxy$ is invariant under $\gamma_\psi\in\Gamtil$,
it is easy to see
that this is independent of the choice of $v$ (\cite[Lemma 2.10]{da}).
In particular,
$$
W_\psi=-\sum_{i=0}^{h_\fp-1}\cxy(e_i).
$$

Let $E$ be an elliptic curve over $F$ of conductor $\mathfrak{p}\nn\infty$ with
split multiplicative reduction at $\mathfrak{p}$ and at $\infty$.
Assume that $E$
is the strong Weil curve in its isogeny class. By the function field
version of the Shimura-Taniyama-Weil conjecture \cite[Section 8]{gr},
there is a newform
$$
c=c_E\in W_{\rsp,S}^{\rm new}(\KK_0(\mathfrak{p}\nn\infty),\ZZ)
$$
which corresponds to $E$. It is uniquely determined up to constants, such that
\begin{enumerate}
\item[(i)]
$c$ is an eigenform under the Hecke algebra with rational eigenvalues and
\item[(ii)]
for all id\`ele class characters $\chi$ of $F$,
$$
L(c,\chi,s)=L(E,\chi,s).
$$
\end{enumerate}
We denote its images in the spaces of harmonic cochains on $\Tp$, $\Tinf$ and
$\Tree$ by $c_\mathfrak{p}$, $c_\infty$ and $c$, respectively.
We choose $c$ such that $c_\infty$ is primitive, i.e., is normalised
such that $c_\infty\in j(\overline{\Gamma_0^\infty(\wp\nn)})$ but 
$c_\infty\not\in nj(\overline{\Gamma_0^\infty(\wp\nn)})$ for $n>1$. This 
determines $c$ up to sign (see section 1.3 and \cite[9.1]{gr}).

\begin{Thm} \label{da_thm}
Put $r={\frac{W_\psi}{\nu_\mathfrak{p}(q_{E,\mathfrak{p}})}}\in\QQ$. Then
$$
I_\psi = \zeta\cdot q_{E,\mathfrak{p}}^r,
$$
where $\zeta$ is a root of unity and $q_{E,\fp}$ the Tate period of
$E$ at $\fp$.
\end{Thm}

\bigskip
Before we can prove theorem \ref{da_thm}, we need to establish a few tools.
A fundamental domain for the action of $\gamma_\psi$ on
$\PP^1(F_\mathfrak{p})-\{\infty,0\}$ is given by
\begin{eqnarray*}
\FF_\psi & = & \bigcup_{i=0}^{h_\fp-1}\pi_\fp^i\OO_\fp^\times \\
& = & \bigcup_{i=0}^{h_\fp-1}U(v_i),
\end{eqnarray*}
where $\pi_\fp$ denotes the matrix
$\begin{pmatrix} \pi_\fp & 0 \\ 0 & 1 \end{pmatrix}$,
and for all $i$,
$$
U(v_i)=\{ t\in\PP^1(F_\mathfrak{p})- \{\infty,0\} : \nu_\mathfrak{p}(t)=i\}.
$$
Each of these sets corresponds to all ends of $\Tp$ which have origin in
$v_i$ and no edge in common with $A(\infty,0)$.

\bigskip
\noindent
We define
$$
m_0:=-\cxy(e_0)=-\cxy(e_{h_\fp})
$$
(since $\cxy$ is invariant under the action of $\gamma_\psi$).

\begin{Prop} \label{I as integral}
For all $z\in\Omega_\mathfrak{p}$,
$$
I_\psi=u^{m_0}\cdot\mint_{\FF_\psi}td\mcxy(t)\in F_\mathfrak{p}^\times.
$$
\end{Prop}

\begin{proof}
The proof proceeds exactly as in \cite[Prop.\ 2.7]{da} (or \cite[Lemma 16]{l}).
\end{proof}

\noindent\textbf{Remark.}
Proposition \ref{I as integral} shows that the definition of $I_\psi$
depends on the choice of $\gamma_\psi$. E.g., if we define a period
$I_\psi'$ using
$$
\gamma_\psi'=\begin{pmatrix} \xi u & 0 \\ 0 & 1 \end{pmatrix},
$$ 
where $\xi$ is a root of unity, we obtain
$I_\psi' = \xi^{m_0}\cdot I_\psi$.

\bigskip
We define $\LL=\{s_1,\dots,s_{h_\infty}\}$ to be the geodesic connecting the
vertices $w_0$ and $w_{h_\infty}$ on $\Tinf$. 
Let $\alpha_1,\dots,\alpha_{h_\infty}\in\Gamtil$ such that
$s_i=\alpha_i s_0$ for $i=1,\dots,h_\infty$.
Furthermore, let 
$$
\lambda:=\sum_{i=1}^{h_\infty}\alpha_i\in\ZZ[\Gamtil]
$$
and
\begin{eqnarray*}
\mu_\LL 
& := & \lambda\ast\mu_{E,\fp} \\
& = &  \sum_{i=1}^{h_\infty}\alpha_i\ast\mu_{E,\fp},
\end{eqnarray*}
where $\mu_{E,\fp}$ is the measure on $\partial\Omega_\fp$
corresponding to $c_\fp$.

\begin{Lemma} \label{measure as sum}
$$
\mcxy=\sum_{n\in\ZZ}\gamma_\psi^n\ast\mu_\LL.
$$
\end{Lemma}

\begin{proof}
Since $A(\infty,0)=\bigcup_{n\in\ZZ}\gamma_\psi^n\LL\subseteq\vec{\EE}(\Tinf)$, 
and using the $\Gamtil$-invariance of $c$, we conclude:
\begin{eqnarray*}
\mcxy(U(e))
& \stackrel{\rm def.}{=} & -\sum_{s\in A(\infty,0)}c(e,s) \\
& = & -\sum_{n\in\ZZ}\sum_{i=1}^{h_\infty}c(e,\gamma_\psi^n\alpha_i s_0) \\
& = & -\sum_{n\in\ZZ}\sum_{i=1}^{h_\infty}c((\gamma_\psi^n\alpha_i)^{-1}e,s_0) \\
& = & \left(\sum_{n\in\ZZ}\gamma_\psi^n\ast\mu_\LL\right)(U(e)).
\end{eqnarray*}
\end{proof}

\begin{Prop} \label{connection with theta functions}
For $i=1,\dots,h_\infty$, let $P_i:=\alpha_i^{-1}(0)$ and
$Q_i:=\alpha_i^{-1}(\infty)$. Then
$$
u^{m_0}\cdot\mint_{\FF_\psi}td\mcxy(t)=
\prod_{i=1}^{h_\infty}\frac{u_{E,\fp}(P_i)}{u_{E,\fp}(Q_i)},
$$
where $u_{E,\fp}$ is the theta function associated to $E$ at $\fp$.
\end{Prop}

\begin{proof}
Similar reasoning as in \cite[Prop.\ 25]{l} shows:
\begin{eqnarray*}
u^{m_0}\mint_{\FF_\psi}td\mcxy(t)
& = & \mint_{\partial\Omega_\fp}td\mu_\LL \\
& \stackrel{\ref{intprop2}}{=} & 
      \prod_{i=1}^{h_\infty}\mint_{\partial\Omega_\fp}\alpha_i(t)d\mu_{E,\fp}(t) \\
& = & \prod_{i=1}^{h_\infty}\mint_{\partial\Omega_\fp}
      \frac{t-P_i}{t-Q_i}d\mu_{E,\fp}(t).
\end{eqnarray*}
The statement follows from Teitelbaum's formula (Prop.\ \ref{poisson}).
\end{proof}

\noindent
Combining propositions \ref{I as integral} and 
\ref{connection with theta functions}, we get
$$
I_\psi=
\prod_{i=1}^{h_\infty}\frac{u_{E,\fp}(P_i)}{u_{E,\fp}(Q_i)}.
$$
Since $P_i$ and $Q_i$ are cusps for all $i$, Theorem \ref{mandri} 
shows that the fraction on the right
is (up to a root of unity) a rational power of the Tate period.
Theorem \ref{da_thm} now follows from the next lemma.

\begin{Lemma}
$\nu_\mathfrak{p}(I_\psi)=W_\psi$.
\end{Lemma}

\begin{proof}
Since $\nu_\fp(u)=h_\fp$,
\begin{eqnarray*}
\nu_\mathfrak{p}(I_\psi) & \stackrel{\ref{I as integral}}{=} & 
      h_\fp m_0 + \nu_\mathfrak{p}(\mint_{\FF_\psi}td\mcxy(t)) \\
& = & h_\fp m_0 + \sum_{i=0}^{h_\fp-1}\mcxy(\pi_\fp^i\OO_\fp^*) \\
& = & h_\fp\cxy(e_0) + \sum_{i=0}^{h_\fp-1}i(\cxy(e_i) - \cxy(e_{i+1})) \\
& = & W_\psi,
\end{eqnarray*}
where the last equality follows from the $\gamma_\psi$-invariance of $\cxy$.
\end{proof}

\medskip
\textbf{Remark:}
Under the assumption that $E/\QQ$ is unique in its $\QQ$-isogeny
class, Darmon (\cite{da}) conjectures that
$$
I_\psi\in q^\ZZ.
$$
His assumption implies especially that $E(\QQ)$ is torsion free.
Since in our case $I_\psi$ is an $F$-rational torsion point,
theorem \ref{da_thm} implies:

\begin{Cor}[Darmon's conjecture] \label{dacon}
Assume that $E(F)_{\rm tor}=0$. Then
$$
I_\psi\in q_{E,\fp}^\ZZ.
$$ 
\end{Cor}

\subsection{The main result}

Finally, we take a look at Teitelbaum's function
field version of the exceptional zero conjecture (\cite{t3}).
He considers the rational function field $F=\FQ(T)$ and the usual place
at infinity $\infty$ with uniformiser $\pi_\infty=\frac{1}{T}$.
We note that in this case $h_\fp=1$ and $W_\psi=m_0$.
We choose $\pi_\fp$ to be an element of $A_\infty=\FQ[T]$.
Since $A_S=\FQ[T,\pi_\fp^{-1}]$, $u$ can chosen to be $u=\pi_\fp$. 

\medskip
Teitelbaum defines a measure on $\OO_\mathfrak{p}$ by
$$
\teit(a+\pi_\mathfrak{p}^n\OO_\pi):=[\frac{a}{\pi_\mathfrak{p}^n},\infty]\cdot c_\infty,
$$
where 
$$
[\frac{a}{\pi_\mathfrak{p}^n},\infty]\cdot c_\infty
=\sum_{s\in A(\frac{a}{\pi_\mathfrak{p}^n},\infty)}c_\infty(s).
$$
Given the elements
\begin{eqnarray*}
m_\mathfrak{p} & = & \nu_\mathfrak{p}(q_{E,\mathfrak{p}}), \\
\tilde{q} & = & \frac{q_{E,\mathfrak{p}}}{\pi_\mathfrak{p}^{m_\mathfrak{p}}}, \\
q(c_\infty) & = & 
\lim_{n\to\infty}\prod_{\stackrel{a \textrm{ mod }\pi_\mathfrak{p}^n}{(a,\pi_\mathfrak{p})=1}}
a^{[\frac{a}{\pi_\mathfrak{p}^n},\infty]\cdot c_\infty} 
= \mint_{\OO_\mathfrak{p}^\times} td\teit(t),
\end{eqnarray*}
he formulates:

\begin{Conj}[Teitelbaum]
$$
\tilde{q}^{[0,\infty]\cdot c_\infty} = q(c_\infty)^{m_\mathfrak{p}}.
$$
\end{Conj}

Our result towards this conjecture is the following:

\begin{Thm} \label{t_conj}
There exists a root of unity $\xi\in \mathbb{F}_{q^{\deg(\fp)}}^\times$ 
such that
$$
\xi\cdot\tilde{q}^{[0,\infty]\cdot c_\infty}
=
q(c_\infty)^{m_\mathfrak{p}}.
$$
\end{Thm}

\noindent
\textbf{Proof of theorem \ref{t_conj}.}
By theorem \ref{da_thm} and proposition \ref{I as integral},
\begin{eqnarray*}
\zeta\cdot (q_{E,\fp})^{\frac{W_\psi}{m_\fp}} 
& = & \pi_\fp^{W_\psi}\mint_{\OO_\fp^\times}td\mcxy(t)
\end{eqnarray*}
for some root of unity $\zeta\in\CC_\fp^\times$.
We set $\xi:=\zeta^{m_\fp}$ and raise to the power $m_\fp$:
\begin{eqnarray*}
\xi\cdot q_{E,\fp}^{W_\psi} 
& = & \pi_\fp^{W_\psi\cdot m_\fp}
\left(\mint_{\OO_\fp^\times}td\mcxy(t)\right)^{m_\fp} \in F_\fp^\times.
\end{eqnarray*}
Also,
\begin{eqnarray*}
W_\psi
& = & -\sum_{s\in A(\infty,0)}c_\infty(s) \\
& = & \sum_{s\in A(0,\infty)}c_\infty(s) \\
& = & [0,\infty]\cdot c_\infty\in\ZZ.
\end{eqnarray*}
This implies that 
$\xi\in \mbox{\boldmath$\mu$}(F_\fp^\times)=\mathbb{F}_{q^{\deg(\fp)}}^\times$, 
since $q_{E,\fp}^{W_\psi}\in F_\fp^\times$.
Dividing by $\pi^{W_\psi\cdot m_\fp}$, the theorem follows from:

\begin{Lemma} \label{meas_comp}
$$
\mcxy\mid_{_{\OO_\mathfrak{p}}}
= \teit.
$$
\end{Lemma}

\begin{proof}
Since $A_\infty$ is dense in $\OO_\fp$,
we only have to consider compact open sets of the form
$$
a+\pi_\fp^n\OO_\fp
=U(\begin{pmatrix} \pi_\fp^n & a \\ 0 & 1 \end{pmatrix}e_0),
$$
where $a\in A_\infty$.

By Lemma \ref{g-inv}, $\mcxy$ is invariant under the action of
$
\begin{pmatrix} -1 & 0 \\ 0 & 1 \end{pmatrix}.
$
In particular,
$$
\mcxy(a+\pi_\fp^n\OO_\fp) = \mcxy(-a+\pi_\fp^n\OO_\fp).
$$
Furthermore,
\begin{eqnarray*}
\mcxy(-a+\pi_\fp^n\OO_\fp)
& = & -\sum_{s\in A(\infty,0)}
      c(\begin{pmatrix} \pi_\fp^n & -a \\ 0 & 1 \end{pmatrix}e_0,s) \\
& \stackrel{(*)}{=} & -\sum_{s\in A(\infty,0)}
      c(e_0,\begin{pmatrix} \pi_\fp^{-n} & \frac{a}{\pi_\fp^n} \\
                              0 & 1 \end{pmatrix}s) \\
& = & -\sum_{s\in A(\infty,\frac{a}{\pi_\fp^n})}c_\infty(s) \\
& = & [\frac{a}{\pi_\fp^n},\infty]\cdot c_\infty \\
& = & \teit(a+\pi_\fp^n\OO_\fp).
\end{eqnarray*}
In $(*)$ we used the $\Gamtil$-invariance of $c$ and the fact that
$\begin{pmatrix} \pi_\fp^n & -a \\ 0 & 1 \end{pmatrix}\in\Gamtil$.
\end{proof}

%% file: doubleint.bbl
\begin{thebibliography}{99}

\bibitem[Da]{da} Darmon, H.: Integration on $\HH_p\times\HH$ and arithmetic
applications, Ann.\ of Math.\ 154 (2001), 589--639.
\bibitem[dS]{ds} de Shalit, E.: $p$-adic periods and modular symbols of
elliptic curves of prime conductor, Inv.\ Math.\ 121 (1995), 225--255.
\bibitem[Dr]{dr} Drinfeld, V.G.: Elliptic modules, Math.\ USSR-Sbornik 23
(1976) 561--592.
\bibitem[Gek1]{gek1} Gekeler, E.-U.: Improper Eisenstein series on 
Bruhat-Tits trees, Manus.\ Math.\ 86 (1995), 367--391.
\bibitem[Gek2]{gek} Gekeler, E.-U.: Analytical Construction of Weil Curves
Over Function Fields, J.\ Th.\ Nomb.\ Bordeaux 7 (1995), 27--49.
\bibitem[Gek3]{gek3} Gekeler, E.-U.: A note on the finiteness of certain
cuspidal divisor class groups, Israel J.\ Math.\ 118 (2000), 357--368.
\bibitem[GR]{gr} Gekeler, E.-U., Reversat, M.: Jacobians of Drinfeld modular
curves, J.\ reine angew.\ Math.\ 476 (1996), 27--93.
\bibitem[GS]{gs} Greenberg, R., Stevens, G.: $p$-adic $L$-functions and
$p$-adic periods of modular forms, Invent.\ Math.\ 111 (1993), no.\ 2,
407--447.
\bibitem[H]{h} Hauer, H.: Teitelbaum's exceptional zero conjecture in the
function field setting, PhD-thesis, Nottingham, 2003.
\bibitem[L]{l} Longhi, I.: Non-Archimedean Integration and Elliptic Curves
over Function Fields, J.\ Numb.\ Th.\ 94 (2002), 375--404.
\bibitem[MTT]{mtt} Mazur, B., Tate, J., Teitelbaum, J.: On $p$-adic
analogues of the conjectures of Birch and Swinnerton-Dyer, 
Invent.\ Math.\ 84 (1986), 1--48.
\bibitem[MT]{mt} Mazur, B., Tate, J.: Refined Conjectures of the
"Birch and Swinnerton-Dyer Type", Duke Math.\ Jour.\ 54 (1987),
711--750.
\bibitem[Se]{se} Serre, J.-P.: \textit{Trees}. Springer-Verlag, 1980.
\bibitem[T1]{t2} Teitelbaum, J.: The Poisson kernel for Drinfeld modular
curves, J.\ Amer.\ Math.\ Soc.\ 4 (1991), 491--511.
\bibitem[T2]{t3} Teitelbaum, J.: Modular symbols for $\FQ[{\rm T}]$,
Duke Math.\ Jour.\ 68 (1992), 271--295.

\end{thebibliography}
